\documentclass[12pt]{article}
\usepackage[cp866]{inputenc}
\usepackage[russian]{babel}

\date{}

\newcounter{theorem}
\newenvironment{theorem}[1][\hspace{-1.0ex}]%
 {\par\addvspace{2mm}\indent\refstepcounter{theorem} \textbf{Теорема\hspace{1.0ex}{\rm#1}.~}\sl}%
 {\par\addvspace{2mm}\rm}
\newcounter{lemma}0
\newenvironment{lemma}[1][\hspace{-1.0ex}]%
 {\par\addvspace{2mm}\indent\refstepcounter{lemma} \textbf{Лемма~\thelemma\hspace{1.0ex}{\rm#1}.~}\sl}%
 {\par\addvspace{2mm}\rm}
\newcounter{corollary}
\newenvironment{corollary}[1][\hspace{-1.0ex}]%
 {\par\addvspace{2mm}\indent\refstepcounter{corollary} \textbf{Следствие~\thecorollary\hspace{1.0ex}{\rm#1}.~}\sl}%
 {\par\addvspace{2mm}\rm}
\newcounter{note}
\newenvironment{note}[1][\hspace{-1.0ex}]%
 {\par\addvspace{2mm}\indent\refstepcounter{note} \textbf{Замечание\hspace{1.0ex}{\rm#1}.~}}%
 {\par\addvspace{2mm}\rm}
\newcounter{example}
\newenvironment{example}[1][\hspace{-1.0ex}]%
 {\par\addvspace{2mm}\indent\refstepcounter{example} \textbf{Пример\hspace{1.0ex}{\rm#1}.~}}%
 {\par\addvspace{2mm}\rm}
\newcounter{problem}
\newenvironment{problem}[1][\hspace{-1.0ex}]%
 {\par\addvspace{2mm}\indent\refstepcounter{problem} \textbf{Проблема~\theproblem\hspace{1.0ex}{\rm#1}.~}}%
 {\par\addvspace{2mm}\rm}

\newcommand\sectionn[1]{\addvspace{5mm}\par\indent{\Large\bf#1}\par\addvspace{2mm}}
\newcommand\sections[1]{\addvspace{5mm}\par\indent\refstepcounter{section}{\Large\bf\thesection.~#1}\par\addvspace{2mm}}

\author{Ю. Л. Васильев, С. В. Августинович, Д. С. Кротов}

\begin{document}

УДК 519.72

\begin{center}
{\LARGE О подвижных множествах\\ в двоичном гиперкубе%
\footnote{~Исследование второго автора выполнено при финансовой поддержке
Российского фонда фундаментальных исследований (проект 07-01-00248)}
}
\end{center}

\begin{center}
{\large Ю. Л. Васильев, С. В. Августинович, Д. С. Кротов}
\end{center}

\begin{abstract}
Если два кода с расстоянием три имеют одинаковую окрестность,
каждый из них называется подвижным множеством.
В двоичном $(4k+3)$-мерном гиперкубе существует подвижное множество мощности
$2\cdot 6^k$, которое нельзя разбить на подвижные множества меньшей мощности
или представить в виде естественного расширения подвижного множества
в гиперкубе меньшей размерности.

If two distance-$3$ codes have the same neighborhood, then each of them is called a mobile set.
In the $(4k+3)$-dimensional binary hypercube, there exists a mobile set of cardinality $2\cdot 6^k$
that cannot be split into mobile sets of smaller cardinalities or represented as a natural extension
of a mobile set in a hypercube of smaller dimension.
\end{abstract}
\sectionn{Введение}

Через $E^n$ обозначается метрическое пространство
на множестве всех двоичных слов
длины $n$ с метрикой Хемминга. 
Пространство $E^n$ иногда называют двоичным, или единичным, или булевым кубом.
Базисный вектор с единицей в $i$-й координате и нулями в остальных обозначается через $e_i$.
Подмножество $M\subseteq E^n$ будем называть $1$-кодом, если
шары радиуса $1$ с центрами из $M$ не пересекаются между собой.
Окрестностью $\Omega(M)$ множества $M$ назовём объединение шаров радиуса $1$ с центрами из $M$,
т.\,е.
$$ \Omega(M) = \{x \in E^n : d(x,M)\leq 1 \}. $$

Если $1$-код $M$ обладает свойством $\Omega(M)=E^n$,
он называется совершенным, или \emph{$1$-со\-вер\-шен\-ным кодом}.
$1$-Совершенные коды существуют лишь в размерностях вида $n=2^k-1$.
Для $n=7$ такой код единственный (с точностью до изометрий пространства)
--- линейный код Хемминга.
При $n=15$ проблема описания и перечисления $1$-совершенных
кодов до сих пор не решена,
несмотря на постоянно растущие возможности вычислительной техники
(существенное продвижение получено в работах \cite{ZZ:2006:l16r14ru,Mal2006:15-15ru}).
В контексте упомянутой проблемы представляется актуальным изучение объектов,
обобщающих в разных смыслах понятие $1$-совершенного кода,
и существующих не только при $n=2^k-1$, но и в промежуточных размерностях.
Такими объектами являются совершенные раскраски
(в частности, двухцветные \cite{FDF:PerfColRU}),
центрированные функции \cite{Vas:2000ACCT}, а также подвижные множества, о которых
пойдёт речь в данной статье.

Множество $M\subseteq E^n$ называется \emph{подвижным} (\emph{п.\,м.}), если:\\
1) $M$ является $1$-кодом,\\
2) существует непересекающийся с $M$ $1$-код $M'$
с той же окрестностью, т.\,е. $M \cap M'=\emptyset$ и $\Omega(M)=\Omega(M')$;
такое множество $M'$ будем называть \emph{альтернативой} множества $M$.

Другими словами, $1$-код есть п.\,м., если у него есть альтернатива.

Для всякого нечётного $n=2m+1$ несложно построить линейное
(замкнутое относительно покоординатного сложения по модулю $2$)
подвижное множество в $E^n$:
\begin{equation}\label{e:lin}
M=\{(x,x,|x|):x\in E^m\}.
\end{equation}
Здесь и далее $|x|$ есть сумма координат вектора $x$ по модулю $2$.
Соответственно
$$ M'=\{(x,x,|x|\oplus 1):x\in E^m\}.$$
Убедиться в выполнении условий $1$ и $2$ для $M$ и $M'$ вполне нетрудно.

Основной целью нашей работы является доказательство следующего факта:

\begin{theorem} {Для всех $n \geq 7$, сравнимых с $3$ по
модулю $4$, в $E^n$ существует нередуцируемое неделимое подвижное
множество.}
\end{theorem}

Непустое п.\,м. $M$ называется \emph{разделимым} (\emph{неделимым}),
если его можно (соответственно, нельзя)
представить в виде объединения двух непустых п.\,м. Понятие
редуцируемости, которое будет сформулировано в разделе~\ref{s:red},
отражает естественную сводимость подвижных множеств к подвижным множествам в размерности на $2$
меньшей.

Простой способ построения п.\,м. в гиперкубах кодовой размерности $n=2^k-1$
заключается в следующем.
Пусть $C$ и $C'$ являются $1$-со\-вер\-шен\-ны\-ми кодами в $E^n$. Тогда $M=C \setminus C'$
есть п.\,м. Действительно, в качестве $M'$ можно взять $C' \setminus C$.
Мощность такого п.\,м. равна $C-|C \cap C'|$.
Мы же исследуем вопрос существования п.\,м.,  которые не сводятся к кодовым размерностям.

В разделе~\ref{s:ext} мы определяем расширенные подвижные множества, в терминах которых
удобно описывать конструкцию. В разделе~\ref{s:i} описана связь подвижных множеств и
$i$-компонент, которые активно изучались ранее.
В разделе~\ref{s:red} описывается конструкция увеличения размерности подвижного множества,
приводящяя к естественному понятию редуцируемого п.\,м.
В разделе~\ref{s:prof} приводится основная конструкция и доказательство теоремы.
В заключительном разделе мы формулируем несколько задач.

\sections{Расширенные подвижные множества} \label{s:ext}

С подвижными множествами, как и с $1$-совершенными кодами, часто бывает удобно работать,
расширив их в следующую размерность проверкой на чётность.
При этом в некоторых случаях получаются более симметричные объекты,
что упрощает доказательства и формулировки утверждений.
И, хотя геометрическая интерпретация расширенных объектов может казаться
не столь изящной и естественной, как в оригинале, и переход с ней требует некоторого
привыкания, многие утверждения становятся более простыми и наглядными,
будучи сформулированными для расширенного случая.

Напомним, что
\emph{расширением} множества $M\subseteq E^n$  называется множество
$\overline M \subseteq E^{n+1}$, полученное добавлением проверки на чётность (нечётность)
ко всем словам множества $M$:
$$ \overline M = \{(x,|x|) : x \in M )\} \qquad\mbox{или}\qquad
\overline M = \{(x,|x|\oplus 1) : x \in M )\} . $$
\emph{Выкалывание} $i$-й координаты в некотором множестве слов из $E^n$
означает удаление $i$-го символа во всех словах множества
(результат будет лежать в $E^{n-1}$). Очевидно, что расширение  и затем выкалывание
последней координаты приводит к исходному множеству,
то есть эти операции в определённом смысле обратные друг другу.

Множество $\overline M\subseteq E^n$ назовём \emph{расширенным подвижным} (\emph{р.\,п.\,м.}),
если оно получается расширением некоторого п.\,м.

Нам будет полезной следующая лемма, которая даёт альтернативные определения р.\,п.\,м.
Как и обычное п.\,м., р.\,п.\,м. $M$ можно определить в паре с другим р.\,п.\,м. $M'$,
которое также естественно называть \emph{альтернативой} р.\,п.\,м. $M$ (из контекста обычно ясно,
речь идёт о подвижных множествах или расширенных подвижных множествах). Для формулировки леммы и дальнейшего использования
удобно определить понятие \emph{сферической окрестности}
$$ \Omega^*(M) =  \Omega(M) \setminus M, $$
которое для расширенных подвижных множеств выполняет
роль, аналогичную роли обычной (``шаровой'') окрестности для п.\,м.
В частности, условие (c) леммы~\ref{l:RPM} определяет р.\,п.\,м. и его альтернативу
аналогично случаю с подвижным множеством.

\begin{lemma}[(об альтернативных определениях р.\,п.\,м.)]\label{l:RPM}
Пусть $M$ и $M'$ есть непересекающиеся $1$-коды в $E^n$,
векторы которых имеют одинаковую чётность $($либо все чётновесовые, либо
нечётновесовые$)$. Пусть $i\in \{1,\ldots,n\}$.
Следующие условия эквивалентны и влекут тот факт, что $M$ $($как и $M')$ есть р.\,п.\,м.

{\rm (a)} Множества $M_i$ и $M'_i$, полученные из $M$ и $M'$ выкалыванием $i$-й координаты,
подвижные
и являются альтернативой друг друга.

{\rm (b)} Граф $($двудольный$)$ расстояний $2$ $G(M \cup M')$
объединения $M \cup M'$ имеет степень $n/2$.

{\rm (c)} $\Omega^*(M) = \Omega^*(M')$.
\end{lemma}

\emph{Доказательство}.
 При $i=n$ из (a) следует, что $M$ есть р.\,п.\,м., по определению. Поскольку условия (b) и (c)
 не зависят от выбора $i$, достаточно показать эквивалентность (a), (b) и (c).

 (c) $\Rightarrow$ (b). Рассмотрим вектор $v$ из $M$.
Рассмотрим множество пар $i,j\in\{1,\ldots,n\}$ таких, что
\begin{equation}\label{e:ij}
v \oplus e_k \oplus e_j \in M'.
\end{equation}
Поскольку $M$ и $M'$ не пересекаются, $k$ и $j$ в такой паре всегда
различны. Поскольку $M'$ есть $1$-код, две различные пары не
пересекаются.
И из условия $\Omega^*(M) = \Omega^*(M')$ следует, что для
любой элемент из $\{1,\ldots,n\}$ принадлежит некоторой паре.
Таким образом, мы имеем разбиение $\{1,\ldots,n\}$ на пары $k,j$,
удовлетворяющие (\ref{e:ij}). Отсюда следует, что степень вершины
$v$ в графе $G(M \cup M')$ равна $n/2$. То же верно для любого $v'$
из $M'$.

 (a) $\Rightarrow$ (c).
 Рассмотрим вектор $w$ на расстоянии $1$ от $M$.
 Нам нужно показать, что он расположен на расстоянии $1$ от $M'$.
 Действительно, в противном случае расстояние от $w$ до $M'$ как минимум $3$
 (учитывая одинаковую чётность $M$ и $M'$),
 и после выкалывания $i$-й координаты он не попадёт в $\Omega(M'_i)$, что противоречит
 условию $\Omega(M_i)=\Omega(M'_i)$.
 Таким образом, $\Omega^*(M)  \subseteq \Omega^*(M')$;
 аналогично, $\Omega^*(M') \subseteq \Omega^*(M) \setminus M$.

 (b) $\Rightarrow$ (a). Рассмотрим вектор $v$ из $M$.
Поскольку степень $v$ в $G(M \cup M')$ равна $n/2$ и в $M'$ нет двух векторов на расстоянии $2$,
все координаты делятся на пары $k$, $j$, удовлетворяющие \ref{e:ij}.
Отсюда, любой вектор вида $v+e_j$, $1\leq j \leq n$, лежит в $\Omega(M')$,
а после выкалывания $i$-й координаты --- в $\Omega(M'_i)$. Но
все такие вектора после выкалывания составляют $\Omega(M_i)$, откуда имеем $\Omega(M_i) \subseteq \Omega(M'_i)$.
Аналогично, $\Omega(M'_i) \subseteq \Omega(M_i)$.
Осталось заметить, что
$M_i \cap M'_i = \emptyset$, поскольку $M$ и $M'$ непересекаются и имеют одну чётность.
$\bigtriangleup$

Учитывая условие (b) и существование линейного п.\,м., имеем следующее важное следствие.

\begin{corollary}
Необходимым и достаточным условием существования непустых п.\,м. $($р.\,п.\,м.$)$ в $E^n$
является нечётность $($соответственно, чётность$)$ $n$.
\end{corollary}

\sections{$i$-Компоненты} \label{s:i}

Содержание данного раздела не используется
при доказательстве основного результата. Однако, оно необходимо для понимания связей
с предшествующими исследованиями, ориентированными на частный случай п.\,м.,
так называемые $i$-компоненты.

П.\,м. $M$ будем называть \emph{$i$-компонентой}, если $\Omega(M)=\Omega(M\oplus e_i)$.
Рассмотрим множество $M_i$, полученное из $M$  выкалыванием $i$-й координаты.
Построим на $M_i$, как на вершинах, так называемый граф минимальных расстояний, соединив
ребром вершины на расстоянии $2$. Доказательство следующей леммы аналогично лемме~\ref{l:RPM},
и мы опускаем его.

\begin{lemma}
\label{l:iK}
\emph{$1$-Код $M$ является $i$-компонентой тогда и
только тогда, когда граф $G(M_i)$ является однородным степени
$(n-1)/2$ и двудольным.}
\end{lemma}

Таким образом, леммы~\ref{l:RPM} и~\ref{l:iK} устанавливают соответствие
между парами альтернативных п.\,м. в $E^{n-1}$
и $i$-ком\-по\-нен\-та\-ми в $E^{n+1}$
(при фиксированном $i$, например, $n+1$).
Это связь проявляется в том, что обоим объектам соответствует
множество в $E^{n}$, граф расстояний два которого является
двудольным и имеет степень $n/2$. В первом случае все вершины
множества будут иметь одинаковую чётность. Во втором --- не
обязательно, однако множества разной чётности будут соответствовать
разбиению $i$-компоненты на независимые $i$-компоненты, ``$i$-чётную'' и ``$i$-нечётную''.
Формально, мы можем сформулировать следующее.

\begin{corollary}\label{c:1}
Множества $M, M' \subseteq E^{n-1}$ есть п.\,м. и его альтернатива если и только если
множество
$$\{(x,|x|,0) : x\in M\} \cup \{(x,|x|,1) : x\in M'\}$$ есть $i$-компонента при $i=n+1$.
\end{corollary}

\begin{corollary}\label{c:2}
Множество $M \subseteq E^{n+1}$ есть $i$-компонента, $i=n+1$, если и только если
множества $$M_a^b = \{x : (x,|x|\oplus a,b)\in M\}, \quad a,b\in\{0,1\}$$
есть п.\,м., причём $M_a^0$ и $M_a^1$ являются альтернативами друг другу
{\rm (множества $M_0^0$ и $M_0^1$ соответствуют ``$i$-чётной'' части $i$-компоненты,
$M_1^0$ и $M_1^1$ --- ``$i$-нечётной''; каждая из этих частей может быть пустой, причём,
если обе непусты, то $i$-компонента разделима).}
\end{corollary}

Примером $i$-компоненты является линейное п.\,м. (\ref{e:lin}), $i=n$.
Ранее \cite{Sol88,Sol:2001i-comp} уже строились многочисленные примеры
нелинейных подвижных множеств, являющихся  $i$-ком\-по\-нентами.
Все они были вложимы в $1$-совершенные коды,
каждый из упомянутых примеров имел мощность, кратную мощности линейной компоненты.
Кроме того, доказанной была лишь неделимость этих $i$-компонент на меньшие
$i$-компоненты. Вопрос об их неделимости как подвижных множеств
остаётся открытым. Поэтому, несмотря на то что исследования посвящены
общей проблеме и даже общему подходу к решению этих проблем,
направления несколько различны и результаты не перекрываются, а дополняют друг друга:
мы отказываемся от вложимости в $1$-совершенные коды (что является ослаблением),
зато имеем дело с более сильной
неделимостью и с б\'{о}льшим спектром размерностей.

\sections{Редуцируемость} \label{s:red}

\begin{lemma}[$($о линейном расширении п.\,м.$)$]
\label{l:ext} Пусть $M,M'\subseteq E^n$ есть р.\,п.\,м. и его альтернатива.
Тогда множества
\begin{equation}\label{e:ext}
R  = \{(x,0,0) : x\in M\} \cup \{(x,1,1) : x\in M'\}
\end{equation}
$$
R' = \{(x,1,1) : x\in M\} \cup \{(x,0,0) : x\in M'\}
$$
есть р.\,п.\,м. и его альтернатива в $E^{n+2}$.
\end{lemma}
\emph{Доказательство}. Выполнение условия (b) леммы~\ref{l:RPM} для
$M$ и $M'$ непосредственно влечёт справедливость этого условия для
$R$ и $R'$. $\bigtriangleup$.

Р.\,п.\,м. $R\in E^n$ назовём \emph{редуцируемым}, если оно может быть
получено конструкцией (\ref{e:ext}), а также перестановкой координат
и инверсией некоторых символов, применёнными ко всем векторам
множества одновременно. П.\,м. назовём \emph{редуцируемым}, если
соответствующее ему р.\,п.\,м. редуцируемо.

Таким образом,   вопрос существования редуцируемых п.\,м. сводится к
существованию п.\,м. в меньших размерностях.
С этой точки зрения формулировка основной теоремы естественна.

\begin{note} Как видно из Следствия~\ref{c:2}, любая $i$-компонента либо
является редуцируемым п.\,м., либо разбивается на две $i$-компоненты
(``$i$-чётную'' и ``$i$-нечётную''), каждая из которых есть редуцируемое п.\,м.
В частности, линейное п.\,м. (\ref{e:lin}) редуцируемое. Более того,
линейное р.\,п.\,м., с точностью до перестановки координат, может быть получено
из тривиального р.\,п.\,м. $\{00\}$ в $E^2$ последовательным применением конструкции
из леммы~\ref{l:ext}.
\end{note}


\sections{Доказательство теоремы} \label{s:prof}

Зафиксируем $n$,  кратное четырём: $n=4k$. Разобьём номера координат на $k$ групп по $4$ в каждой
и переобозначим соответствующие орты следующим образом:
$e_0^1, e_1^1, e_2^1, e_3^1, e_0^2, \ldots , e_3^k$.
В каждой четвёрке вида
$\{e_0^i, e_1^i, e_2^i, e_3^i\}$ выберем произвольно (всего $6$ возможностей)
 пару несовпадающих ортов $e_j^i$ и $e_t^i$ и назовём её индексом число $p=j \star t-1$,
 где $\star$ определяется таблицей значений
 $$
 \begin{array}{c|cccc}
 \star & 0 & 1 & 2 & 3 \cr
 \hline
 0 & 0 & 1 & 2 & 3 \cr
 1 & 1 & 0 & 3 & 2 \cr
 2 & 2 & 3 & 0 & 1 \cr
 3 & 3 & 2 & 1 & 0 \cr
 \end{array}
 $$
 (заметим, что $j \star t = t \star j$ и $a \star b = c \star d$ для любых попарно различных
 $a$, $b$, $c$, $d$).
 Просуммировав выбранные пары по всем $i=1,2,\ldots,k$, мы получим вектор веса $2k$,
 который будем называть \emph{стандартным}. Всего получится $6^k$ стандартных векторов.
 Индексом $I(V)$ стандартного вектора $v$ мы будем называть сумму по модулю $3$ всех индексов
 составляющих его пар ортов.

 Разобьём множество стандартных векторов на непересекающиеся подмножества $S_0$, $S_1$ и  $S_2$
 в соответствии с их индексами.

 \textbf{Утверждение~1.} \emph{Пусть $i\neq j$, $i,j \in \{0,1,2\}$.
 Тогда граф $G(S_i \cup S_j)$ расстояний два,
 индуцированный множеством векторов $S_i \cup S_j$, является двудольным и однородным степени $2k$.}

 Для начала заметим, что графы $G(S_i)$ и $G(S_j)$ пусты. Действительно, рассмотрим пару векторов
 $v,u \in S_i$.  Либо $v$ и $u$ различаются в одной четвёрке координат,  тогда $d(v,u)=4$,
 поскольку у них индексы одинаковы, либо $v$ и $u$ различаются в большем, чем одна, числе четвёрок,
 тогда $d(v,u) \geq 4$, поскольку по каждой четвёрке расстояние между стандартными векторами чётно.
 Таким образом, двудольность графа $G(S_i \cup S_j)$ обеспечена.

 Также легко понять, что всякий вектор индекса $i$ имеет ровно двух
 соседей на расстоянии $2$ из $S_j$, различающихся с ним в одной
 фиксированной четвёрке координат. Это означает, что степень графа
 есть $2k$. Утверждение~1 доказано.

 Таким образом, множество $S_0$ (например) является
 расширенным подвижным и имеет мощность $2 \cdot 6^{k-1}$.

 \textbf{Утверждение~2.} \emph{Р.\,п.\,м. $S_0$ неделимое.}

 Предположим,  что $P \subseteq S_0$ и $S_0 \backslash P$ есть р.\,п.\,м. и $P$ непусто.
 Тогда существует $P$ имеет альтернативу $P'$.

 Сначала убедимся, что\\
 (*) \emph{$P'$ состоит из стандартных векторов}, то есть
 таких, что в каждой четвёрке координат 
 содержится ровно две единицы. 
 Действительно, в противном случае
 $P'$ содержит вектор с нестандартной четвёркой,
 и, как следствие, $\Omega^*(P')$ содержит вектор с двумя нестандартными
 четвёрками.
 В то же время $\Omega^*(P)$
 состоит из векторов с одной нестандартной четвёркой
  и, следовательно, не может совпадать с
 $\Omega^*(P')$,
 что противоречит лемме~\ref{l:RPM}. (*)
 доказано.

 Рассмотрим произвольный вектор $p$ из $P$ и покажем, что\\
 (**) \emph{все векторы из $S_0$,  отличающиеся от $p$
 не более чем в двух четвёрках, также принадлежат $P$.}
 Без потери общности рассмотрим две первые четвёрки. Положим
 $p=(h,t)$, где $h$ и $t$ --- векторы длины $8$ и $n-8$ соответственно.
 Рассмотрим вектор $p \oplus e_j^i$ из $\Omega^*(P)$, где
 $i\in\{1,2\}$ и $j\in\{0,1,2,3\}$. Согласно леммы~\ref{l:RPM},
 $p \oplus e_j^i \in \Omega^*(p')$ для некоторого $p'$ из $P'$.
 Как доказано выше, вектор $p'$ стандартный, поэтому $p'=p \oplus e_j^i \oplus
 e_{j'}^i$ для некоторого $j' \in \{0,1,2,3 \}$, откуда следует, что
 $p'$ совпадает с $p$ в последних $n-8$ координатах.
 Из этих рассуждений следует, что $\Omega^*(P_8) = \Omega^*(P'_8)$, где
 $$ P_8 = \{ b \in E^8 : (b,t)\in P \} $$
 $$ P'_8 = \{ b \in E^8 : (b,t)\in P' \} $$
 и, по утверждению (c) леммы~\ref{l:RPM}, множество $P_8$ есть р.\,п.\,м. в $E^8$.
 Легко установить (например, пользуясь утверждением (b) леммы~\ref{l:RPM}),
 что мощность р.\,п.\,м. в $E^8$ больше $6$.
 С другой стороны, по построению, ровно $12$ векторов из $S_0$
 имеют вид $(b,t)$, $b \in E^8$. Следовательно, больше половины
 таких векторов принадлежат $P$. Если бы не все принадлежали $P$, то
 к оставшимся векторам (из $S_0 \backslash P$) были бы применимы аналогичные
 рассуждения, что привело бы к противоречию. Следовательно, все $12$
 векторов из $S_0$, совпадающих с $p$ во всех координатах кроме
 первых восьми, принадлежат $P'$, что доказывает (**).

 Таким образом, любые два вектора из $S_0$ на расстоянии $4$ друг от
 друга одновременно либо  принадлежат $P$, либо нет. Поскольку
 $S_0$, очевидно, связно по расстоянию $4$, получаем $P=S_0$.
 Утверждение~2 доказано.

 \textbf{Утверждение~3.} \emph{Р.\,п.\,м. $S_0$ не редуцируемо.}

 Заметим, что в конструкции (\ref{e:ext}) сумма последних двух
 координат равна $0$ для любого слова из $R$. Учитывая перестановку
 координат и инверсию символов, можно утверждать, что у редуцируемого
 р.\,п.\,м. существуют две координаты, сумма которых равна $0$ либо $1$
 одновременно для всех слов множества. Легко проверить, что $S_0$ не
 удовлетворяет этому условию: любые две координаты содержат все
 четыре комбинации из $0$ и $1$. Утверждение~3 доказано. Теорема доказана.

\sections{Заключение}

Мы построили бесконечный класс неделимых нередуцируемых п.\,м.
Конструкция обобщает пример, упомянутый в конце работы \cite{VasSol97ru}.
В заключение мы сформулируем несколько задач,
 естественно связанных с исследованием подвижных множеств и с проблемой характеризации
 их многообразия.

Для построения п.\,м. можно применять принцип обобщённой каскадной
конструкции для $1$-совершенных кодов
\cite{ZinLob:2000ru}. В частности, конструкция из раздела~\ref{s:prof}
может быть интерпретирована в таких терминах.
Неделимые п.\,м., построенные таким образом, будут иметь неполный ранг, то
есть
для всех слов множества координаты будут удовлетворять некоторому линейному уравнению
(неполной проверке на чётность или нечётность).

\begin{problem}
Построить бесконечный класс неделимых п.\,м. полного ранга.
\end{problem}

\begin{example}
Рассмотрим четыре слова
$$
\begin{array}{r@{\,}c@{\,}l}
  (&100 \cr &110 \cr &010&)
\end{array}
,\quad
\begin{array}{r@{\,}c@{\,}l}
  (&011 \cr &110 \cr &000&)
\end{array}
,\quad
\begin{array}{r@{\,}c@{\,}l}
  (&101 \cr &001 \cr &011&)
\end{array}
,\quad
\begin{array}{r@{\,}c@{\,}l}
  (&001 \cr &100 \cr &111&)
\end{array}
,
$$
из $E^9$, записанные для удобства в виде массива $3\times 3$, а также все
слова, полученные из них циклическими перестановками строк и/или
столбцов массива.
Получим неделимое п.\,м. полного ранга мощности $36$.
Альтернатива получается инверсией всех слов.
\end{example}

\begin{problem}
Построить богатый класс транзитивных неделимых п.\,м., р.\,п.\,м.
Множество $M\subseteq E^n$ называется
транзитивным, если стабилизатор ${\mathrm{Stab}}_I(M)$ множества $M$ в группе $I$ изометрий гиперкуба
действует транзитивно на элементах $M$,
т.\,е. для любых $x,y$ из $M$ найдётся изометрия $\sigma\in \mathrm{Stab}_I(M)$
такая, что $\sigma(x)=y$.
Например, нетрудно показать, что п.\,м., построенные в данной работе, являются транзитивными.
Известно несколько конструкций транзитивных $1$-совершенных и расширенных $1$-совершенных кодов,
последние результаты смотри в \cite{Sol:2005transitiveRU,Pot:transRU}.
\end{problem}

\begin{problem}
 Исследовать вопрос вложимости п.\,м. в $1$-совершенный код:
существование невложимых п.\,м. в кодовых размерностях $n=2^k-1$;
существование п.\,м., невложимых при помощи линейного расширения (лемма~\ref{l:ext})
в $1$-совершенный код ни в одной большей размерности.
В частности, для п.\,м., построенных в разделе~\ref{s:prof}, вопросы вложимости открыты при $n\geq 11$.
\end{problem}

\begin{problem}
 Оценить максимальный размер неделимого п.\,м.
\end{problem}

\begin{problem}
 Оценить минимальный размер нелинейного п.\,м. (конструкция раздела~\ref{s:prof} вместе с леммой~\ref{l:ext}
 даёт верхнюю оценку $1{,}5L(n)$, где $L(n)=2^{(n-1)/2}$ --- мощность линейного п.\,м.),
 нередуцируемого неделимого п.\,м. (конструкция даёт верхнюю оценку $1{,}5^{(n-3)/4}L(n)$),
 неделимого п.\,м. полного ранга.
\end{problem}


\renewcommand\refname{\indent Литература}


\end{document}


УДК 519.72

\begin{center}
{\LARGE On mobile sets in the binary hypercube%
\footnote{~The second author is partially supported by the RFBR grant 07-01-00248}
}
\end{center}

\begin{center}
{\large Yu. L. Vasil'ev, S. V. Avgustinovich, D. S. Krotov}
\end{center}

\renewcommand\abstractname{\indent Abstract}
\begin{abstract}
If two distance-$3$ codes have the same neighborhood, then each of them is called a mobile set.
In the $(4k+3)$-dimensional binary hypercube, there exists a mobile set of cardinality $2\cdot 6^k$
that cannot be split into mobile sets of smaller cardinalities or represented as a natural extension
of a mobile set in a hypercube of smaller dimension.
\end{abstract}
\sectionn{Introduction}

By $E^n$ we denote the metric space of all length-$n$ binary words with the Hamming metric.
The space $E^n$ is called the binary, or unary, or Boolean hypercube.
The basis vector with one in the $i$th coordinate and zeros in the other is denoted by $e_i$.
A subset $M$ of $E^n$ is called a $1$-code if
the radius-$1$ balls with centers in $M$ are disjoint.
The union of the radius-$1$ balls with the centers in $M$ is called the \emph{neighborhood} of $M$ and denoted by
$\Omega(M)$, i.e.,
$$ \Omega(M) = \{x \in E^n : d(x,M)\leq 1 \}. $$

If a $1$-code $M$ satisfies $\Omega(M)=E^n$, then it is called perfect, or a \emph{$1$-per\-fect code}.
$1$-Per\-fect codes exist only when the dimension has the form $n=2^k-1$.
For $n=7$, such a code is unique (up to isometries of the space),
the linear Hamming code.
For $n=15$, the problem of characterization and enumeration of the $1$-perfect
codes is not solved yet, in spite of the increasing computation abilities
(considerable results are obtained in \cite{ZZ:2006:l16r14,Mal2006:15-15}).
In this context, it is topical to study objects that generalize, in different senses,
the concept of $1$-perfect code and exist in intermediate dimensions, not only of type $n=2^k-1$.
Examples of such objects are the perfect colorings
(in particular, with two colors \cite{FDF:PerfCol}),
the centered functions \cite{Vas:2000ACCT}, and the mobile sets, discussed in this paper.

A set $M\subseteq E^n$ is called \emph{mobile} (\emph{m.s.}) iff:\\
1) $M$ is a $1$-code;\\
2) there exists  a  $1$-code $M'$  disjoint with $M$ and
with the same neighborhood, i.e., $M \cap M'=\emptyset$ and $\Omega(M)=\Omega(M')$;

such a set $M'$ will be called the \emph{alternative} of $M$.

In other words, a $1$-code is a m.s. iff it has an alternative.

For every odd $n=2m+1$, we can construct a linear
(closed  with respect to coordinatewise modulo-$2$ addition)
m.s. in $E^n$:
\begin{equation}\label{e:lin}
M=\{(x,x,|x|):x\in E^m\}.
\end{equation}
(Here and below $|x|$ denotes the modulo-$2$ sum of the coordinates of $x$.)
Respectively,
$$ M'=\{(x,x,|x|\oplus 1):x\in E^m\}.$$
It is not difficult to check the conditions $1$ and $2$ for these $M$ and $M'$.

Our main goal is to prove the following:

\begin{theorem} {For all $n \geq 7$ congruent to $3$ modulo $4$,
there exists an irreducible unsplittable mobile set in $E^n$.}
\end{theorem}

A nonempty m.s. $M$ is called \emph{splittable} (\emph{unsplittable}),
iff if can (respectively, cannot) be represented as the union of two nonempty m.s.
The concept of reducibility, which will be defined in Section~\ref{s:red},
reflects a natural  reducibility of mobile sets to mobile sets in the hypercube of
the two-less dimension.

A simple way to construct a m.s. in a hypercube of a code dimension $n=2^k-1$
is the following.
Let  $C$ and $C'$ are $1$-perfect codes in $E^n$. Then $M=C \setminus C'$
is a m.s. Indeed, we can take $C' \setminus C$ as $M'$.
The cardinality of this m.s. is $C-|C \cap C'|$.
We study the existence of m.s. that cannot be reduced to code dimensions.

In Section~\ref{s:ext} we define extended mobile sets; that concept is convenient for the description
of our construction.
In Section~\ref{s:i} we describe a connection between the mobile sets and the $i$-components,
which were studied earlier.
In Section~\ref{s:red} we describe a construction of increasing dimension for mobile sets;
that construction leads to the natural concept of a reducible m.s.
In Section~\ref{s:prof} we give the main construction and prove Theorem.
In the final section, we formulate several problems.

\sections{Extended mobile sets} \label{s:ext}

Like as with $1$-perfect codes, it is sometimes convenient to work with mobile sets extending them
by  the all-parity check to the next dimension.
In some cases we get more symmetrical objects,
which simplifies proofs and formulations of statements.
And.
Some statements become more simple and intuitive while being formulated for the extended case,
although geometrical interpretations of extended objects can seem to be not so elegant and natural
as for the original.

Recall that the
\emph{extension} of the set $M\subseteq E^n$  is the set
$\overline M \subseteq E^{n+1}$ obtained by the addition of the all-parity-check bit
to all the words of $M$:
$$ \overline M = \{(x,|x|) : x \in M )\} \qquad\mbox{or}\qquad
\overline M = \{(x,|x|\oplus 1) : x \in M )\} . $$
 \emph{Puncturing}  the $i$th coordinate for some set of words in $E^n$
means  removing the $i$th symbol from all the words of the set
(the result is in $E^{n-1}$).
Obviously, the extension and puncturing the last coordinate lead to the original set;
so, these operations are opposite to each other, in some sence.

A set $\overline M\subseteq E^n$ is called \emph{extended mobile} (an \emph{e.m.s.})
iff it can be obtained as the extension of some m.s.

We will use the following lemma, which gives alternative definitions of an e.m.s.
As a usual m.s., an e.m.s. $M$ can be defined together with some other e.m.s. $M'$,
which can also be referred as an \emph{alternative} of $M$
(usually, it is clear from the context what we are talking about, mobile sets or extended mobile sets).
For the formulation of the lemma and further using, it is convenient to define the concept
of the \emph{spherical neighborhood}
$$ \Omega^*(M) =  \Omega(M) \setminus M, $$
which, for the extended mobile sets, plays the role similar to the role of the usual (``ball'') neighborhood for the m.s.
In particular, part (c) of Lemma~\ref{l:RPM} defines an e.m.s. and an alternative similarly
to the case of a m.s.

\begin{lemma}[(alternative definitions of an e.m.s.)]\label{l:RPM}
Let $M$ and $M'$ be disjoint $1$-codes in $E^n$,
and let their vectors have the same parity $($either all vectors are even, or odd$)$.
Let $i\in \{1,\ldots,n\}$.
The following conditions are equivalent and imply
that $M$ $($as like as $M')$ is an e.m.s.

{\rm (a)} The sets $M_i$ and $M'_i$ obtained from $M$ and $M'$ by puncturing  $i$th coordinate
are mobile and, moreover, are alternatives of each other.

{\rm (b)} The $($bipartite$)$ distance-$2$ graph $G(M \cup M')$
of the union $M \cup M'$ has the degree $n/2$.

{\rm (c)} $\Omega^*(M) = \Omega^*(M')$.
\end{lemma}

\emph{Proof}.
[The translation of the proof to appear]
$\bigtriangleup$

Taking into account (b) and the existence of a linear m.s., we have the following important corollary.

\begin{corollary}
Nonempty m.s. $($e.m.s.$)$ exist in $E^n$ if and only if
$n$ is odd $($respectively, even$)$.
\end{corollary}

\sections{$i$-Components} \label{s:i}


A m.s. $M$ is called an \emph{$i$-component} iff $\Omega(M)=\Omega(M\oplus e_i)$.
Consider the set $M_i$ obtained from $M$ by puncturing the $i$th coordinate.
Построим на $M_i$, как на вершинах, так называемый граф минимальных расстояний, соединив
ребром вершины на расстоянии $2$. The proof of the following lemma is similar to Lemma~\ref{l:RPM},
and we omit it.

\begin{lemma}
\label{l:iK}
\emph{A $1$-code $M$ is an $i$-component if and only if
the graph $G(M_i)$ is regular of degree $(n-1)/2$ and bipartite.}
\end{lemma}

So, Lemmas~\ref{l:RPM} and~\ref{l:iK} establish a correspondence
between pairs of alternative m.s. in $E^{n-1}$
and $i$-com\-po\-nents in $E^{n+1}$
(for fixed $i$, say, $i=n+1$).
This correspondence is evident as both objects correspond to a set in $E^{n}$
whose distance-$2$ graph is bipartite and has the degree  $n/2$.
In the first case, all the vertices of this set have the same parity.
In the second case, this is not not necessary,
but the subsets of different parity will correspond to a partition
of the $i$-component into independent $i$-components, ``$i$-even'' and ``$i$-odd''.
Formally, we can formulate the following.

\begin{corollary}\label{c:1}
Sets $M, M' \subseteq E^{n-1}$ are a m.s. and an alternative if and only if
the set
$$\{(x,|x|,0) : x\in M\} \cup \{(x,|x|,1) : x\in M'\}$$
is an $i$-component with $i=n+1$.
\end{corollary}

\begin{corollary}\label{c:2}
A set $M \subseteq E^{n+1}$ is an $i$-component with $i=n+1$ if and only if
the sets $$M_a^b = \{x : (x,|x|\oplus a,b)\in M\}, \quad a,b\in\{0,1\}$$
are m.s., where $M_a^0$ and $M_a^1$ are alternatives to each other
{\rm (the sets $M_0^0$ and $M_0^1$ correspond to the ``$i$-even'' part of the $i$-component;
$M_1^0$ and $M_1^1$, to the ``$i$-odd''; each of these parts can be empty;
and if both are nonempty, then the $i$-component is splittable).}
\end{corollary}

An example of $i$-component is the linear m.s. (\ref{e:lin}), $i=n$.
Formerly \cite{Sol88en,Sol:2001i-comp} many examples of nonlinear $i$-components were constructed.
Each of them is embeddable to a $1$-perfect codes and has the cardinality,
divisible by the cardinality of the linear component.
Moreover, it was only proved that these $i$-components cannot be split into smaller
$i$-components.
Their splittability onto mobile sets are still questionable.
So, in spite of the fact that the researches are devoted to common problems
and a common approach,
the lines are slightly different and the results do not overlap but  complement each other:
we give the embeddability to $1$-perfect codes up (which is a weakening)
but deal with a stronger splittability and a wider specter of dimensions.

\sections{Reducibility} \label{s:red}

\begin{lemma}[$($on the linear extension of a m.s.$)$]
\label{l:ext} Let $M\subseteq E^n$ be an e.m.s. and let $M'\subseteq E^n$ be an alternative of $M$.
Then the set
\begin{equation}\label{e:ext}
R  = \{(x,0,0) : x\in M\} \cup \{(x,1,1) : x\in M'\}
\end{equation}
is an e.m.s. with an alternative
$$
R' = \{(x,1,1) : x\in M\} \cup \{(x,0,0) : x\in M'\}.
$$
\end{lemma}
\emph{Proof}. Condition (b) of Lemma~\ref{l:RPM} for
$M$ and $M'$ implies the validity of this condition for
$R$ and $R'$. $\bigtriangleup$.

An e.m.s. $R\in E^n$ is called \emph{reducible} iff it can be obtained by the construction (\ref{e:ext})
and applying some isometry of the space (i.e., a coordinate permutation and the inversion in some coordinates).
A m.s. is called \emph{reducible} iff
the corresponding e.m.s. is reducible.

So, the existence of reducible m.s. is reduced to
 the existence of m.s. in smaller dimensions.
From this point of view, the formulation of the main theorem is natural.

\begin{note} As we can see from Corollary~\ref{c:2}, any $i$-component is either
reducible m.s. or can be split into two  $i$-components
(``$i$-even'' and ``$i$-odd''), which are reducible m.s.
In particular, the linear m.s. (\ref{e:lin}) is reducible. Moreover,
the linear e.m.s., up to a coordinate permutation, can be obtained from
the trivial e.m.s. $\{00\}$ in $E^2$ by  sequential applying the construction from
Lemma~\ref{l:ext}.
\end{note}

\sections{Proof of Theorem} \label{s:prof}

Let us fix $n$ divisible by $4$: $n=4k$.
Partition the coordinate numbers into $k$ groups with $4$ numbers in each group;
rename the corresponding orts as follows:
$e_0^1, e_1^1, e_2^1, e_3^1, e_0^2, \ldots , e_3^k$.
In each quadruple  of type
$\{e_0^i, e_1^i, e_2^i, e_3^i\}$ we chose arbitrarily (there exist $6$ pocibilities)
 a pair of different orts $e_j^i$ and $e_t^i$; by the index of the pair we shell mean the number $p=j \star t-1$,
 where $\star$ is defined by the value table
 $$
 \begin{array}{c|cccc}
 \star & 0 & 1 & 2 & 3 \cr
 \hline
 0 & 0 & 1 & 2 & 3 \cr
 1 & 1 & 0 & 3 & 2 \cr
 2 & 2 & 3 & 0 & 1 \cr
 3 & 3 & 2 & 1 & 0 \cr
 \end{array}
 $$
 (note that $j \star t = t \star j$ and $a \star b = c \star d$ for any distinct
 $a$, $b$, $c$, $d$).
 Summarizing the chosen pairs for all $i=1,2,\ldots,k$, we get a vector of weight $2k$,
 which will be called \emph{standard}. Totally, there exist $6^k$ standard vectors.
 By the index  $I(v)$ of a standard vector $v$ we shell mean the modulo-$3$ sum of
 the indexes of all the pairs of orts that constitute $v$.

 Let us partition the set of standard vectors into disjoint subsets $S_0$, $S_1$, and  $S_2$
 in compliance with the indexes of the vectors.

 \textbf{Claim~1.} \emph{Let $i\neq j$, $i,j \in \{0,1,2\}$.
 Then the distance-two graph $G(S_i \cup S_j)$
 induced by the set $S_i \cup S_j$ is bipartite and regular of degree $2k$.}

[The translation of the proof to appear] Claim~1  is proved.

 So, $S_0$ (for example) is an e.m.s. of cardinality $2 \cdot 6^{k-1}$.

 \textbf{Claim~2.} \emph{The e.m.s. $S_0$ is unsplittable.}

 [The translation of the proof to appear]
 Claim~2  is proved.

 \textbf{Claim~3.} \emph{The e.m.s. $S_0$ is irreducible.}

[The translation of the proof to appear] Claim~3 is proved. The theorem is proved.

\sections{Conclusion}

We have constructed an infinite class of unsplittable irreducible m.s.
Our construction generalizes the example mentioned in \cite{VasSol97}.
In conclusion, we formulate several problems, which are
 naturally connected with the study of mobile sets and with the problem of characterization of their
variety.

For constructing m.s., one can apply the generalized concatenation principle
which works for $1$-perfect codes
\cite{ZinLob:2000}. In particular, the construction from Section~\ref{s:prof}
can be treated in such terms.
Unsplittable m.s. constructed in such the way will have non-full rank,
i.e., for all the words of the set the coordinates will satisfy some linear equation.

\begin{problem}
Construct an infinite family of full-rank unsplittable m.s.
\end{problem}

\begin{example}
Consider the four words
$$
\begin{array}{r@{\,}c@{\,}l}
  (&100 \cr &110 \cr &010&)
\end{array}
,\quad
\begin{array}{r@{\,}c@{\,}l}
  (&011 \cr &110 \cr &000&)
\end{array}
,\quad
\begin{array}{r@{\,}c@{\,}l}
  (&101 \cr &001 \cr &011&)
\end{array}
,\quad
\begin{array}{r@{\,}c@{\,}l}
  (&001 \cr &100 \cr &111&)
\end{array}
,
$$
in $E^9$, listed, for convenience, as $3\times 3$ arrays, and all the
words obtained from them by cyclic permutations of rows and/or
columns of the array.
We get full-rank unsplittable m.s. of cardinality $36$.
An alternative can be obtained by the inversion of all the words.
\end{example}

\begin{problem}
Construct a rich class of transitive unsplittable m.s., e.m.s.
A set $M\subseteq E^n$ is called
transitive iff the stabilizer ${\mathrm{Stab}}_I(M)$ of $M$ in the group $I$ of isometries
of the hypercube acts transitively on the elements of $M$;
i.e., for every $x,y$ from $M$ there exists an isometry $\sigma\in \mathrm{Stab}_I(M)$
such that $\sigma(x)=y$.
For example, it is not difficult to see that the m.s. constructed in the current paper are transitive.
There are several constructions of transitive $1$-perfect and extended $1$-perfect codes,
see \cite{Sol:2005transitive,Pot:trans} for the last results.
\end{problem}

\begin{problem}
 Study the embeddability of m.s. into $1$-perfect codes:
the existence of nonembeddable m.s. in the code dimensions $n=2^k-1$;
the existence of m.s. that cannot be embedded with help of the linear extension (Lemma~\ref{l:ext})
into a $1$-perfect code in a larger dimension.
In particular, for m.s. constructed in Section~\ref{s:prof}, the embedding questions are open provided $n\geq 11$.
\end{problem}

\begin{problem}
 Estimate the maximal  cardinality of an unsplittable m.s.
\end{problem}

\begin{problem}
 Estimate the minimal cardinality of a nonlinear m.s.
 (the construction of Section~\ref{s:prof} together with Lemma~\ref{l:ext}
 give the upper bound $1{,}5L(n)$, where $L(n)=2^{(n-1)/2}$ is the cardinality of the linear m.s.),
 of an irreducible unsplittable m.s. (the construction gives the upper bound $1{,}5^{(n-3)/4}L(n)$),
 unsplittable m.s. of full rank.
\end{problem}


\renewcommand\refname{\indent References}


\providecommand\href[2]{#2} \providecommand\url[1]{\href{#1}{#1}}